\documentclass[a4paper,draft,oneside,12pt]{article}
\usepackage{ amsthm,amsmath,amsfonts,amssymb}
\usepackage{epsfig}
\usepackage{colordvi,helvet,multicol}
\usepackage{epic}
\usepackage{mathptmx}
\usepackage{times}
\usepackage[T1]{fontenc}
\usepackage{verbatim}
\usepackage{graphics}
\usepackage{makeidx}
\usepackage{fancyhdr}
\makeindex
\makeatletter
\long\def\@caption#1[#2]#3{%        Abbildungstext auf Fussnotengroesse
  \par
  \addcontentsline{\csname ext@#1\endcsname}{#1}%
    {\protect\numberline{\csname the#1\endcsname}{\ignorespaces #2}}%
  \begingroup
    \@parboxrestore
    \if@minipage
      \@setminipage
    \fi
   \normalsize
    \footnotesize
    \@makecaption{\csname fnum@#1\endcsname}{\ignorespaces #3}\par
  \endgroup}
\makeatother
 \theoremstyle{plain}
\newcommand{\ee}{\end{equation}}

\newtheorem{definition}{Definition}[section]
\newtheorem{theorem}[definition]{Theorem}

\newtheorem{remark}[definition]{Remark}

\newtheorem{corollary}[definition]{Corollary}

\newtheorem{pr1}[definition]{}
\newcommand\RR{{\Bbb R}}

\newcommand\NN{{\Bbb N}}

\begin{document}
\title{ Hypergroups with  Unique $\alpha$-Means 
}

 \author{Ahmadreza Azimifard \\
\footnotesize\texttt{}
}

\date{}
\maketitle

\begin{description}
    \item Abstract:
    
                 Let $K$ be a commutative hypergroup and $\alpha\in \widehat{K}$. We
                 show that  $K$ is $\alpha$-amenable with the unique $\alpha$-mean
                 $m_\alpha$ if and only if $m_\alpha\in L^1(K)\cap L^2(K)$ and 
                 $\alpha$ is isolated in $\widehat{K}$. In contrast to the case of 
                 amenable noncompact locally compact   groups,  examples of
                 polynomial  hypergroups with unique $\alpha$-means ($\alpha\not=1$)
                 are given.  Further examples   emphasize
                 that  the $\alpha$-amenability of hypergroups depends
                 heavily on the asymptotic behavior of  Haar measures   and  characters.
                 \footnote{
Parts of the paper
    are taken from     the author's   Ph.D.~thesis at the Technical University of Munich.}

     \item R\'esum\'e:
     
	              Soit $K$  un hypergroupe commutatif et  $\alpha\in \widehat{K}$ .
     						Nous montrons que   $K$ est $\alpha$-moyennable avec unicité de l'$\alpha$-moyenne
    						$m_\alpha$ si et seulement si $m_\alpha\in L^1(K)\cap L^2(K)$ et
    						$\alpha$ est isolé dans $\widehat{K}$. 
    						Contrairement au cas des groupes moyennables localement compacts mais non compacts,
    						des exemples d'hyper-groupes polynomiaux avec unicité des $\alpha$-moyennes
    						($\alpha\not=1$) sont donnés. Nous montrons à l'aide d'autres examples 
    						que l'$\alpha$-moyennabilité des hypergroupes dépend 
    						fortement de leurs mesures de Haar ainsi que du comportement 
    						du caractères.

\end{description}

{
\footnotesize{
										\begin{tabular}{lrl}
										\bf  {Keywords.}  &  \multicolumn{2}{l}
										 { Hypergroups:  orthogonal polynomial,  of Nevaei classes.}\\
										 &{ $\alpha$-Amenable Hypergroups.}\\
										\end{tabular}
										\vspace{.2cm}


}}

%     \end{pr1}

 %%%%%%%%%%%%%%%%%%%%%%%%%%%%%%%%%%%%%%%%%%%%%%%%%%%%%%%%%%%%%%%%%%%%%%%%%%%%%%%%%%%%%%%%%

\newpage
 							
%%%%%%%%%%%%%%%%%%%%%%%%%%%%%%%%%%%%%%%%%%%%%%%%%%%%%%%%%%%%%%%%%%%%%%%%%%%%%%%%%%%%%%
	\section{Introduction}						
					Recently the notion of $\alpha$-amenable hypergroups was 
								 introduced and studied in \cite{f.l.s}. Let $K$ be a commutative locally compact  hypergroup and 
								let $L^1(K)$ 
								 denote  the hypergroup algebra.
								 Assume that  
								 $\alpha\in \widehat{K}$ and denote by $I(\alpha)$ 
								 the maximal ideal in  $L^1(K)$ generated by $\alpha$.
    As shown in \cite{f.l.s},
 								$K$ is $\alpha$-amenable if and only if either  
 								$I(\alpha)$ has a b.a.i. (bounded approximate identity) or  $K$ satisfies
 								the modified Reiter's condition of $P_1$-type in  $\alpha$.
 								  Commutative hypergroups are always 1-amenable \cite{Ska92}, whereas 
 								  a large class of
 								  non $\alpha$-amenable 
											 hypergroups, $\alpha\not=1$, are given in \cite{thesis, f.l.s}. 
											It is worth to notify   that $1\in \mbox{ supp }\pi_K$ does not 
											 hold in general, where 
											 $\mbox{ supp }\pi_K$ denotes  the support
											 of the Plancherel measure on $\widehat{K}$ \cite{Jew75, Ska92}.

											  As in the case of   locally compact groups \cite{Pat88},
								 if $K$ is a  noncompact locally compact amenable  hypergroup, then 
                 the cardinality of (1-)means is $2^{2^{d}}$, 
                 where $d$ is the smallest cardinality 
                 of a cover of $K$ by compact sets \cite{Ska92}.
				However, it is well known that    $K$ has a unique (1-)mean if and only 
				         if $K$ is compact \cite{Pat88, Ska92}. Hence, 
								 $\mbox{ supp }\pi_K=\widehat{K}$ and $K$ is 
				 $\alpha$-amenable for every $\alpha\in \widehat{K}$ \cite{BloHey94, f.l.s}.

 							 For a $\alpha$-amenable hypergroup $K$ with a unique $\alpha$-mean,
 							 one can pose the natural question of whether  $K$
 							 is compact or  $K$ is $\beta$-amenable when $\alpha\not=\beta \in \widehat{K}$.
 				Theorem \ref{main.theorem} answers
 				this question completely.  In addition, examples
 				of polynomial hypergroups show that 
 				the $\alpha$-amenability of hypergroups depends  on
 				the asymptotic behavior of the Haar measures and 
 				characters. Furthermore, the $\alpha$-amenability of  $K$ with
 				a unique $\alpha$-mean ($\alpha\not=1$), even   in every $\alpha\in \widehat{K}\setminus{\{1\}}$,
 				does not imply the compactness of $K$;  see Section \ref{examples}.		
							
	Different axioms for  hypergroups  are given 
							in \cite{ Dun73, Jew75, Spec75}. However, in this paper 
							we refer to Jewett's axioms in \cite{Jew75}.
%%%%%%%%%%%%%%%%%%%%%%%%%%%%%%%%%%%%%%%%%%%%%%%%%%%%%%%%%%%%%

																			\section{Preliminaries}
					       Let $(K, \omega, \sim )$ be a locally compact
					       hypergroup, where 
                          $\omega:K\times K\rightarrow M^1(K)$ defined by $(x,y)\mapsto \omega(x,y)$,
                and $\sim:K\rightarrow K$ defined by $x\mapsto \tilde{x}$, 
                denote the convolution and involution on $K$, 
                where  $M^1(K)$ stands for the set of  all probability 
                measures on $K$. The hypergroup $K$ is called commutative 
                if  $\omega{(x,y)}=\omega{(y,x)}$  for every $x, y\in K$. 
     				
     			Throughout this  paper $K$ is a commutative hypergroup.
     			Let $C_c(K)$ be the	spaces of all  continuous 
     			functions on $K$ with compact support.
          The translation of $f\in C_c(K)$ at  
    			the point $x\in K$, $T_xf$,  is defined by
      \begin{equation}\notag
      						T_xf(y):=\int_K  f(t)d\omega{(x,y)}(t), \mbox{ for every } y\in K.
      \end{equation}
        	Being $K$ commutative  ensures the existence of a
        	Haar measure $m$ on $K$ which is unique up to a 
        	multiplicative constant \cite{Spec75}.
         Let         $(L^p(K), \|\cdot \|_p)$ $(p=1, 2)$  
         denote the usual Banach space  of Borel measurable functions on $K$ \cite[6.2]{Jew75}.
        For   $f, g\in L^1(K)$  we may define  the convolution and
        involution    by
                              $ f*g(x):=\int_K f(y)T_{\tilde{y}}g(x)dm(y)$ ($m$-a.e. on $K$)
                                    and
                         $f^\ast(x)=\overline{f(\tilde x)}$, 
                          respectively, that     $\left(L^1(K), \|\cdot \|_1\right)$ becomes a   Banach $\ast$-algebra.
                      %%%%%%%%%%%%%%%%%%%%%%%%%%%%%%%%%%%%%%%%%%%%%%%%%%%%%%%%%%%%%
         %
         If $K$ is discrete, then $L^1(K)$ has an 
         identity element. Otherwise $L^1(K)$ has a
         b.a.i.,  i.e. there exists a net $\{e_i\}_i$ 
         of functions in $L^1(K)$ with $\|e_i\|_1\leq M$,
         for some $M>0$, such that $\|f \ast e_i-f\|_1\rightarrow 0$ 
         as $i\rightarrow \infty$  \cite{BloHey94}.
                         The set of all multiplicative linear functionals on 
                         $L^1(K)$, i.e.  the maximal ideal space of $L^1(K)$  \cite{BonDun73},
                         can be identified  with
                   \begin{equation}\notag
                           \mathfrak{X}^b(K):=\left\{
                                                    \alpha\in C^b(K): \alpha\not=0, \;\omega(x,y)
                                                    (\alpha)=\alpha(x)\alpha(y), \; \forall\;x,y\in K
                                             \right\}
                  \end{equation}
                                   via
                                   $ \varphi_\alpha(f):=\int_K f(x)\overline{\alpha(x)}dm(x)$, for every  $f\in L^1(K)$.
                                   $\mathfrak{X}^b(K)$  is a 
                                   locally compact Hausdorff space with  
                                    the compact-open topology \cite{BloHey94}.
                                    $\mathfrak{X}^b(K)$ and 
                                    its subset
                                    \[
                                     \widehat{K}:=\{\alpha\in \mathfrak{X}^b(K):
                                     \alpha(\tilde x)=\overline{\alpha(x)}, \;  \forall x\in K\}
                                   \]
                    are considered as the character spaces of $K$.
                    The maximal ideal in $L^1(K)$ generated by the character $\alpha$  is
                    $I(\alpha):=\{f\in L^1(K):\varphi_\alpha(f) =0\}$.
                    The Fourier transform of $f\in L^1(K)$, 
                    $\widehat{f}\in C_0(\widehat{K})$, is
                    $\widehat{f}(\alpha):=\varphi_\alpha(f)$ 
                    for every   $\alpha\in \widehat{K}$.
        There exists a unique (up to a multiplicative constant) 
        regular  positive Borel measure $\pi_K$ on $\widehat{K}$ 
        with $\mbox{ supp }\pi_K=\mathcal{S}$ 
        such that 
        $\int_K |f(x)|^2dm(x)=\int_{\mathcal{S}}|\widehat{f}(\alpha)|^2 d\pi_K(\alpha)$ 
        for all $f\in L^1(K)\cap L^2(K) $  \cite{BloHey94}.
         The   extension of 
         the Fourier transform 
         defined on $L^1(K)\cap L^2(K)$ 
         to all of  $L^2(K)$ onto $L^2(\widehat{K})$ is the Plancherel transform 
         which  is an isometric isomorphism.
         Observe that $\mathcal{S}$ is a 
         nonvoid closed subset of $\widehat{K}$,
         and  the constant function 1 is 
         in general not contained in $\mathcal{S}$ \cite[9.5]{Jew75}.

        The inverse Fourier transform for 
        $\varphi\in L^1(\widehat{K})$ is given  by
        $\check{\varphi}(x)=\int_{\mathcal{S}} \varphi(\alpha)\alpha(x)d\pi_K(\alpha)$
        for every  $x\in K$. Then 
        $\check{\varphi}\in C_0(K)$ and 
        if $\check{\varphi}\in L^1(K)$ 
        then 
        $\widehat{\check{\varphi}}=\varphi$ \cite{BloHey94}.

					Let $L^1(K)^\ast$ and $L^1(K)^{\ast\ast}$ denote the dual 
					and the bidual spaces  of $L^1(K)$ respectively. As usual, 
					$L^1(K)^\ast$ can be  identified with the space $L^\infty(K)$ of essentially bounded Borel 
					measurable complex-valued functions on $K$.
 					We may define     the Arens product  on $L^1(K)^{\ast\ast}$  as follows:
											\begin{equation}\notag
														\langle m\cdot m',  f\rangle=\langle m, m'\cdot f\rangle
											\end{equation}
             in which
						$\langle m'\cdot f,  g \rangle=\langle m',  f\cdot g\rangle$
						and $\langle f\cdot g,  h\rangle=\langle f,  g \ast h\rangle$
						for all $m, m'\in L^1(K)^{\ast\ast}$, $f\in L^\infty(K)$ 
						and $g,h\in L^1(K)$. 
			       $L^1(K)^{\ast\ast}$ with the Arens product is 
			       a noncommutative Banach algebra in general \cite{BonDun73, Civ}.
        From   the definitions of the Arens product and the convolution 
        we may have 
        $g\cdot f={g}^{\ast}\ast f$ and   $m\cdot(f \cdot g)=(m\cdot f)\cdot g.$
        
%%%%%%%%%%%%%%%%%%%%%%%%%%%%%%%%%%%%%%%%%%%%%%%%%%%%%%%%%%%%%%%%%%%%%%%%%%%%%%%%%%
  
 \begin{definition}\label{ch.2.14}
           \emph{Let $K$ be a commutative hypergroup 
                 and $\alpha \in \widehat{K}$. $K$ is 
                 called $\alpha$-amenable if  there 
                 exists a bounded linear functional 
                 $m_\alpha$ on $L^{\infty }({{K}})$ with
                 the following properties:}
\begin{itemize}
              \item[\emph{(i)}]$m_\alpha(\alpha)=1$,
                  \item[\emph{(ii)}]$m_\alpha({\delta_{\tilde x}}\ast f)
                 							       =
                          						{\alpha(x)} m_\alpha (f), $ \hspace{.2in}\emph{for
                          							every }$f\in L^{\infty }({{K}})$ \emph{and }$x \in
                          							K$.
                          \end{itemize}
 \end{definition}

            For example, if $K$ is compact or $L^1(K)$ is amenable, then $K$ 
             is $\alpha$-amenable, for every $\alpha\in \widehat{K}$ \cite{f.l.s, Ska92}.

 \section{Main Theorem}

				\begin{theorem}\label{main.theorem}
												\emph{
															Let $K$ be a hypergroup and $\alpha\in \widehat{K}$.
															If $K$ is 
															$\alpha$-amenable with the  
															unique $\alpha$-mean $m_\alpha$, then}
									\begin{itemize}
													\item[\emph{(i)}] \emph{
																									$m_\alpha$ and $\alpha$ belong to 
																									$L^1(K)\cap L^2(K)$ and 
																									$\alpha\in \mathcal{S} $ 
																									 is isolated. Further,
																									$m_\alpha^2=m_\alpha$.
																									}
													\item[\emph{(ii)}] \emph{
																										$m_\alpha=\pi(\alpha)/\|\alpha\|_2^2$, 
																										where 
																										$\pi:L^1(K)\rightarrow L^1(K)^{\ast\ast}$ 
																										is the canonical embedding.
																									}
													\item[\emph{(iii)}] \emph{
																												If $\alpha$ is positive,  
																												 then $\alpha=1$, hence $K$ is compact.
																									 }
                  \end{itemize}
      \end{theorem}
%______________________________________________________________________________
 
\begin{proof}
								Since $K$ is $\alpha$-amenable with 
								the unique $\alpha$-mean $m_\alpha$,
								$m_\alpha(\alpha)=1$, and $ f\cdot{g} ={g}\cdot f$, we have 
								\begin{align}\notag
                   \langle m_\alpha, f\cdot g \rangle  &=\langle m_\alpha, {g}^\ast \ast f\rangle\\ \notag 
                                                    &=\langle m_\alpha, \int_K (\delta_{x}\ast f) g^\ast(x) dm(x)\rangle \\ \notag 
                                                   & =\int_K \langle m_\alpha,\delta_{x}\ast f\rangle g^\ast(x) dm(x) \\\notag
                                                   & =\widehat{g^{\ast} }(\alpha) \langle m_\alpha, f\rangle, \notag
               \end{align}
						for every $f\in L^{\infty}(K)$ and $g\in L^1(K)$. 
						Moreover, if  $n\in L^1(K)^{\ast \ast}$ and $ h\in L^1(K)$, 
						then
          \[
            \langle m_\alpha\cdot n,f\cdot{g} \rangle 
              =
                \langle m_\alpha ,n\cdot(f\cdot{g})\rangle
              =
                \langle m_\alpha, (n\cdot f)\cdot{g}\rangle
              = 
                {\widehat{g^{\ast}} (\alpha)}
                \langle m_\alpha, n\cdot f\rangle
              =
                {\widehat{g^\ast} (\alpha)} 
                \langle m_\alpha\cdot n ,f\rangle.
          \]

										Since the $\alpha$-mean $m_\alpha$ is unique and 
										the associated functional to $\alpha$ 
										on $L^1(K)^{\ast\ast}$ 
										is multiplicative \cite{Civ}, 
										$m_\alpha\cdot n=\lambda_n\cdot m_\alpha$,
										where $\lambda_n=\langle n, \alpha\rangle$.
										Let $(n_i)$ be a net in $L^1(K)^{\ast\ast}$
										converging to $n$ in the $w^\ast$-topology.
										Then the convergence   $\lambda_{n_i}\rightarrow \lambda_n$, as $i\rightarrow \infty$,  implies that
										the mapping  $n\rightarrow {m_\alpha}\cdot n$   is 
										$w^*$-$w^*$ continuous on 
										$L^1(K)^{\ast\ast}$,
										hence
									$m_\alpha$ is in  $L^1(K)$,  the topological centre of 
									$L^1(K)^{\ast\ast}$ \cite{Kamyabi.2}.
 				In that $\widehat{m_\alpha}(\alpha)=1$, 
  			$g\cdot m_\alpha={\widehat{g^\ast}(\alpha)}m_\alpha$
  			for every $g\in L^1(K)$, 
  			and the Arens product is continuous
  			in the first variable, then  
  			$m_\alpha^2=m_\alpha.$

 											Let $\beta\in \widehat{K}$. 
 											The equality
														$\beta(x)m_\alpha(\beta)
														=m_\alpha(T_x \beta)
														=\alpha(x)m_\alpha(\beta)$,
														 for all $x \in K$,
										implies that    
									 $m_\alpha(\beta)=\delta_\alpha(\beta)$.
 														Since  $\widehat{m_\alpha}\in C_0(\widehat{K})$,  $\alpha$ is isolated in
														$\widehat{K}$ and $\widehat{m_\alpha}\in L^1(\widehat{K})$.
													  The inverse of Fourier theorem yields  
													  $m_\alpha=\widehat{m_\alpha}{^\vee}$, hence
													  $\alpha\in \mathcal{S}$. Moreover, 
													  since  the Plancherel  transform is an  isometric 
													  isomorphism  of  $L^2(K)$
													  onto  $L^2(\widehat{K})$
													  and $\widehat{m_\alpha}(\beta)=\delta_{\beta}(\alpha)$, 
													    $m_\alpha\in L^2(K)$.\\
 
(ii)						 Plainly  $\delta_x\cdot m_\alpha=\alpha(x)m_\alpha$, 
								for every $x\in K$,  so  it follows from part (i) that  
								 $\alpha\in L^1(K)\cap L^2(K)$.

									Let $n_\alpha=\pi(\alpha)/\|\alpha\|_2^2$.
									We shall  prove  $m_\alpha=n_\alpha$. 
									Apparently  $\langle n_\alpha,\alpha \rangle=1$, 
									and for every 
									$x\in K$ and $f\in L^\infty(K)$ 
									we have $\langle n_\alpha, T_x f\rangle =\alpha(x)\langle n_\alpha, f\rangle $,
									hence $n_\alpha$ is a $\alpha$-mean on 
									$L^\infty(K)$. 
									Since   $m_\alpha\in L^1(K)^{\ast\ast}$, there exists 
									$(m_i)_i$ a net of
									functions in $L^1(K)$
									such that $\pi({m_j})\overset{w^\ast}{\longrightarrow} m_\alpha$, 
								  Goldstein's theorem  \cite{Dunford}. Moreover, 
									 ${m_j}\cdot {\alpha}=\alpha \cdot m_j
									 ={\widehat{{m_j}}(\alpha)}{{\alpha}}$ and 
									 $m_\alpha(\alpha)=1$, 
									 so  taking the $w^\ast$-limit    yields 
									 $m_\alpha\cdot  {\pi(\alpha)}=\pi(\alpha)$. Therefore,
									  for every $f\in L^\infty(K)$ and $x\in K$
									  we have
										\[
										  \|{\alpha}\|_2^2  \langle  n_{{\alpha}}, 
										  f\rangle= \langle {\pi({\alpha})} , f \rangle
										 =
										  \langle   m_\alpha\cdot\pi({\alpha}),f\rangle
								    =
								    \langle m_\alpha,\pi({\alpha})\cdot f \rangle
								    =
								    \langle m_\alpha,{\alpha}\cdot f \rangle 
								    =
								    \|{\alpha}\|_2^2 \langle m_\alpha, f \rangle,
  									\]
   hence
$m_\alpha=n_{{\alpha}}$.

(iii) 										By (i) since $\alpha\in L^1(K)\cap L^2(K)$, 
													we have
														 \[
														 		\alpha(x)\int_K \alpha(y)dm(y)
														 		=
														 		\int_KT_x\alpha(y)dm(y)=\int_K \alpha(y)dm(y),
														 \]
													which implies that $\alpha(x)=1$ for 
													every $x\in K$,  hence $K$ is compact \cite{Jew75}.

\end{proof}

%%%%%%%%%%%%%%%%%%%%%%%%%%%%%%%%%%%%%%%%%%%%%%%%%%%%%%%%%%%%%%%%%%%%%%%%%%

					\begin{corollary}
												\emph{
																Let  $K$ be a  $\alpha$-amenable hypergroup 
																with a unique $\alpha$-mean in all	
																$\alpha\in \widehat{K}\setminus\{1\}$. 
																Then $1\in \mathcal{S}$.
																}
					\end{corollary}

		\begin{remark}
									\emph{ 
									We observe that part (iii) of Theorem \ref{main.theorem} can also be derived from part (i) and 
									\cite[Theorem 2.1]{voit.p.c.}.
												 											}
   \end{remark}

        \section{ Examples}\label{examples}
        
        \begin{itemize}

           \item[(I)]{\bf{
                          Symmetric hypergroup \cite{Voit91.1}:
												 }
																}
																	For each $n\in \NN$,  let $b_n\in ]0,1]$, $c_0=1$, 
																	and 
																	define numbers $c_n$ inductively 
																	by $c_n=\frac{1}{b_n}\left(c_0+c_1+...c_{n-1}\right)$.
 													A symmetric 
 													hypergroup structure on $\NN_0$ is 
 													defined by 
 													$\varepsilon_n\ast \varepsilon_m=\varepsilon_m\ast \varepsilon_n=\varepsilon_n  \mbox{ if  } 0\leq m<n$ 
 													and 
 												\[\varepsilon_n\ast \varepsilon_n
 													  =
 													  \frac{c_0}{c_n}\varepsilon_0+\frac{c_1}{c_n}\varepsilon_1
 													  +
 													  ...\frac{c_{n-1}}{c_n}\varepsilon_n+(1-b_n)\varepsilon_n.\]
 													  
 													$\NN_0$ with the above  convolution and 
 													an  involution defined by the identity map    
 													is a commutative  hypergroup with  
 													$\mathfrak{X}^b(\NN_0)=\mathcal{S}$.
 													Every  nontrivial character $\alpha$ in 
 													$ \widehat{\NN_0}$ has 
 													a finite  support, so 
 													$\alpha\in \ell^1(\NN_0)\cap \ell^2(\NN_0)$. 
 													Consequently   by Theorem \ref{main.theorem} we see that  
 													 $\NN_0$ is  $\alpha$-amenable with a unique $\alpha$-mean if and only if $\alpha\not=1$.
        
       \item [(II)]                           
                            Let $\{p_n\}_{n\in \NN_0}$ be a   set  of   polynomials defined by a recursion relation
     \begin{equation}\label{r.10}
                       p_1(x)p_n(x)=a_np_{n+1}(x)+b_np_n(x)+ c_np_{n-1}(x)       
		\end{equation}
								 for $n\in \NN$ and $p_0(x)=1$, $p_1(x)=\frac{1}{a_0}(x-b_0)$, 	where 
									$a_n>0$,  $b_n\in \RR$ for all $n\in \NN_0$ and $c_n > 0$ 
 									for 
 									$n\in \NN$. There exists a probability measure 
 									$\pi\in M^1(\RR)$ such that 
 									$\int_\RR p_n(x)p_m(x)d\pi(x)
 									=
 									\delta_{n,m}\mu_m\hspace{.1in}(\mu_m>0)$ \cite{Chi78}.
 									 Assume that  $p_n(1)\not=0$, so  after renorming, 
 									  for $n\in \NN_0$ we have $p_n(1)=1$. 
 									  The relation (\ref{r.10}) implies that  
 									  $a_n+b_n+c_n=1$ and $a_0+b_0=1$. 
 									  The  polynomial set   $\{p_n\}_{n\in \NN_0}$   
 									   induces a hypergroup structure 
 									   on $\NN_0$ \cite{BloHey94}, 
 									   which is known as a polynomial hypergroup.

    \begin{itemize}
           \item [(i)] {\bf{
 																		 Hypergroups of compact type \cite{Frankcompact}:
 																	  }
 															  }
																 If in  the    recursion 
																			formula (\ref{r.10})  $a_n, c_n\rightarrow 0$ and $b_n\rightarrow 1$ 
																			    as $n\rightarrow \infty$, then the induced  hypergroup $\NN_0$ is called to be of compact type.
																			    In this case,   $\mathcal{S}=\widehat{\NN_0}=\mathfrak{X}^b(\NN_0)$, 1 is the 
																			    only accumulation point of $\widehat{\NN_0}$ and nontrivial characters of $\NN_0$ belong to 
																			    $\ell^1(\NN_0)\cap \ell^2(\NN_0)$. By Theorem \ref{main.theorem}
																			    we see that $\NN_0$ is $\alpha$-amenable with a unique $\alpha$-mean if and only if $\alpha\not=1$.
																			     For instance, the  little q-Legendre polynomial    hypergroup is  of compact type.
										 
  						 \item[(ii)]{
  						 								\bf{
  						 											Hypergroups of Nevai Classes:
  						 										}
  						 										}\label{examples.nevai}
															Let $\{p_n\}_{n\in \NN_0}$  define  a hyergroup structure on $\NN_0$ with the 
															relations  (\ref{r.10}).
															Consider  the orthonormal polynomials  
															$q_n(x):=\sqrt{h(n)}p_n(x)$, which  by the 
															recursion (\ref{r.10}) satisfy the following 
															recursion formula 
								\begin{equation}\notag
								               xq_n(x)=\lambda_{n+1}q_{n+1}(x)
								               +
								               \beta_nq_n(x)+\lambda_n q_{n-1}(x),\hspace{1.cm} \forall n\in \NN_0,
							\end{equation}
							 where $q_0(x)=1$, $\lambda_n=a_0\sqrt{c_na_{n-1}}$ for $n\geq 2$,
							  $\lambda_1=a_0\sqrt{c_1}$, 
							  $\lambda_0=0$, and $\beta_n=a_0b_n+b_0$ for 
							  $n\geq 1$, with $\beta_0=b_0$. 
							  The polynomial set $(q_n)_{n\in \NN_0}$ is of
							  the Nevai class $M(0,1)$ if 
							  $\underset{n\rightarrow \infty}{\lim}\lambda_n=\frac{1}{2}$
							   and
							   $\underset{n\rightarrow \infty}{\lim}\beta_n=0$. 
							   It has a  bounded variation, $(q_n)_{n\in \NN_0}\in BV$,
							   if 
							   \[\sum_{n=1}^{\infty}
							   (|\lambda_{n+1}-\lambda_n|+|\beta_{n+1}-\beta_n|)<\infty.\]

               %  Observe that if    $(q_n)_{n\in \NN_0} \in M(0,1)$, then  
               %       $\mathcal{S}=[-1,1]\cup A$, 
               %  where $A$ is a nonempty  countable set 
               %  and $[-1, 1]\cap A= \emptyset$, see  \cite[Theorem 7]{Nev79}.

              \begin{theorem}\label{main.2}
                            \emph{
                                   Let $(q_n)_{n\in \NN_0}\in BV\cap M(0,1)$ and $\alpha_x\in\widehat{\NN_0}$, where  
                                   $\alpha_x(n):=q_n(x)$ for $n\in \NN_0$. Then the followings hold:
                        \begin{enumerate}
                               \item[(i)] $\mathcal{S} \cong [-1,1]\cup A$, 
                													 where $A$ is a nonempty  countable set 
                														 and $[-1, 1]\cap A= \emptyset$.
                               \item[(ii)] 
                                            If $x\in A$, then  $\NN_0$ is $\alpha_x$-amenable 
                                            with a unique $\alpha_x$-mean.
                               \item[(iii)]  If $h(n)$ is unbounded, 
                                            then $\NN_0$ is not $\alpha_x$-amenable for  $x\in (-1,1)$.
                               \item[(iv)] If $h(n)$ is bounded, 
                                            then $\NN_0$ is $\alpha_x$-amenable for  $x\in (-1,1)$.
                        \end{enumerate}
                                 }
              \end{theorem}
%________________________________________________________________________________________________

     \begin{proof}
     							(i) It is shown in \cite[Theorem 7]{Nev79}.
     							
									(ii)  
									Let $\mathcal{S}$ be as in  part (i). 
									If
									 $A\cap ]1,\infty [\not=\emptyset$, then $x_1:=\mbox{sup}A\in A$
									  corresponds to a positive
									  character of $\NN_0$, \cite[Theorem 5.3]{Chi78}.
									  But this contradicts the fact that
									  a positive character in  $\mathcal{S}$ 
									  cannot be isolated, \cite[Theorem 2.1]{voit.p.c.},
									  hence 
									  $A\subset ]-\infty, -1[$. 
									  By \cite[Theorem 18(p.36)]{Nev79}  we have
                   \[
                     \underset{n\rightarrow \infty}{\lim}\frac{h(n+1)}{h(n)}
                           \left|
                                \frac{p_{n+1}(x)}{p_{n}(x)}
                          \right|
                                =
                                 C\underset{n\rightarrow \infty}{\lim}\left|
                                                                            \frac{p_{n+1}(x)}{p_{n}(x)}
                                                                     \right|
                               = 
                                   \left(
                                         |x|+(x^2-1)^{1/2}
                                  \right)^{-1}
                                  <1,
                   \]
                       whenever $x\in A$. 
                       This  shows that $\alpha_x$  belongs 
                       to $\ell^1(\NN_0)$. 
                       Hence,   by Theorem \ref{main.theorem}  
                       $\NN_0$ is $\alpha_x$-amenable with  a  unique $\alpha_x$-mean.

                      (iii) and (iv) are   shown in \cite[Theorems 4.10--11]{f.l.s}.
        \end{proof}
 %%%%%%%%%%%%%%%%%%%%%%%%%%%%%%%%%%%%%%%%%%%%%%%%%%%%%%%%%%%%%%%%%%%%%%%%%%%%%%%%%
 
							\begin{remark}
                    \begin{enumerate}
                            \item[\emph{(i)}] \emph{
                                                  Theorem \ref{main.2} reveals    that the 
                                                  $\alpha$-amenability of 
                                                  $K$ in general  depends  on the asymptotic 
                                                  behavior of the Haar measure  and  
                                                  $\alpha$.
                                                   }
                            \item[\emph{(ii)}]\emph{  
                                                     Observe that in     
                                                     Theorem \ref{main.2} (iii) if
                                                      $x\in (-1,1)$ then   the functionals $m_{\alpha_x}$ are
                                                      distinct.
                                                   }
                   \end{enumerate}
            \end{remark}
   \end{itemize}
   \end{itemize}
%%%%%%%%%%%%%%%%%%%%%%%%%%%%%%%%%%%%%%%%%%%%%%%%%%%%%%%%%%%%%%%%%%%%%%

       \begin{pr1}{\bf{Conjecture:}}
                       \emph{
                             Let $K$ be a $\alpha$-amenable hypergroup. Then $K$ has either a unique $\alpha$-mean
                              or 
                              the cardinality of the set 
                             of $\alpha$-means is at most  $2^{2^d}$, 
                             where $d$ is the smallest 
                             cardinality of  a cover  of $K$ by compact sets.
                             }

      \end{pr1}
      
%%%%%%%%%%%%%%%%%%%%%%%%%%%%%%%%%%%%%%%%%%%%%%%%%%%%%%%%%%%%%%%

%%%%%%%%%%%%%%%%%%%%%%%%%%%%%%%%%%%%%%%%%%%%%%%%%%%%%%%%%%%%%%%%%%%%%%%%%%%%%%%%%%%%%%%%%%%%%%%%%%%%%%%%%%%%%%%

%E-Mail: azimifard@yahoo.com
 
\end{document}